\documentclass[12 pt]{article}
\usepackage{amscd, amssymb,amsmath,amsthm}
\newtheorem{lm}{Lemma}
\newtheorem{tm}{Theorem}

\newtheorem{Corollary}{Corollary}

\newcommand{\ext}{\operatorname{Ext}}

\newcommand{\T}{\mathbf{T}}
\newcommand{\A}{\mathbf{A}}

\newcommand{\proj}{\mathbf{P}}
\newcommand{\Z}{\mathbb{Z}}

\newcommand{\cO}{\mathcal{O}}
\newcommand{\cI}{\mathcal{I}}

\newcommand{\bV}{\mathsf{V}}

\newcommand{\bw}{\mathsf{w}}
\newcommand{\bW}{\mathsf{W}}
\newcommand{\bv}{\mathsf{v}}

\newcommand{\bF}{\mathsf{F}}
\newcommand{\bZ}{\mathsf{Z}}

\newcommand{\rarr}{\rightarrow}
\newcommand{\oh}{{\mathcal{O}}}
\newcommand{\com}{\mathbb{C}}
\newcommand{\Q}{\mathbb{Q}}

\newcommand{\eqq}{\stackrel{\sim}{=}}

\newcommand{\bpf}{\noindent {\em Proof.} }
\newcommand{\epf}{\qed \vspace{+10pt}}

\newcommand{\ch}{\text{ch}}

\begin{document}

\title{Gromov-Witten theory and Donaldson-Thomas theory, II}
\author{D.~Maulik, N.~Nekrasov, A.~Okounkov, and R.~Pandharipande}
\date{1 July 2004}

\maketitle

%\vspace*{-8.5cm}
%\hfill{\scriptsize{ITEP-TH-61/03}} 
%
%\hfill{\scriptsize{\ IHES/M/03/67}} 
%\vspace*{8.5cm}

\pagestyle{plain}
\section{Introduction}

\subsection{Overview}

The Gromov-Witten theory of a 3-fold $X$ is defined via integrals over the 
moduli space of stable maps. The Donaldson-Thomas theory of $X$
is defined via integrals over the moduli space of ideal sheaves.
In \cite{MNOP1}, a GW/DT correspondence equating the two theories was proposed,
and the Calabi-Yau case was presented. 
We discuss here the GW/DT correspondence for
general 3-folds.

Let $X$ be a nonsingular, projective 3-fold. Insertions in the
Gromov-Witten theory of $X$ are determined by primary and descendent fields. 
Insertions in the 
Donaldson-Thomas theory of $X$ are naturally obtained
from the Chern classes of universal
sheaves. We conjecture a GW/DT correspondence for 3-folds
relating these two sets of insertions. 

Let $S\subset X$ be a nonsingular surface. The Gromov-Witten theory of
$X$ relative to $S$ has been defined in \cite{EGH,IP,LR,L}. 
The relative constraints
are determined by partitions weighted by cohomology classes of $S$.
A relative Donaldson-Thomas theory has been defined by J. Li \cite{jli}.
The relative constraints are determined by cohomology classes of the
Hilbert scheme of points of $S$.
We propose a GW/DT correspondence in the relative case relating the
Gromov-Witten constraints to the Donaldson-Thomas constraints via Nakajima's 
basis of the cohomology of the Hilbert scheme of points.

In the last Section of the paper, independent of the conjectural framework,
we study the Donaldson-Thomas theory in degree 0 using localization
and relative geometry. We derive a formula for the equivariant vertex measure 
in the degree 0 case and prove Conjecture $1'$
of \cite{MNOP1} in the toric case. A degree 0 relative formula is also proven.

\subsection{Acknowledgments}

We thank J. Li for explaining his definition of relative Donaldson-Thomas
theory to us. An outline of his ideas is presented in Section \ref{relid}.
We thank J. Bryan, T. Graber,
A. Iqbal, M. Kontsevich, Y. Soibelman, R. Thomas, and C. Vafa for 
related discussions.

D.~M. was partially supported by a Princeton Centennial graduate fellowship.
N. ~N. was partially supported 
by the grants RFFI 03-02-17554 and NSh-1999.2003.2.
He is grateful to the Princeton Mathematics department
for hospitality.
A.~O.\ was partially supported by
DMS-0096246 and fellowships from the Sloan and Packard foundations.
R.~P.\ was partially supported by DMS-0071473
and fellowships from the Sloan and Packard foundations.

\section{The GW/DT correspondence for 3-folds}

\subsection{GW theory}
\label{gwt}
Gromov-Witten theory is defined via integration over the moduli
space of stable maps.
Let $X$ be a nonsingular, projective 3-fold. Let
 $\overline{M}_{g,r}(X,\beta)$ denote the moduli space of
$r$-pointed stable maps from connected, genus $g$ curves to $X$ representing the
class $\beta\in H_2(X, \Z)$. Let 
$$\text{ev}_i: \overline{M}_{g,r}(X,\beta) \rarr X,$$
$$ L_i \rarr \overline{M}_{g,r}(X,\beta)$$
denote the evaluation maps and cotangent lines bundles associated to
the marked points.
Let $\gamma_1, \ldots, \gamma_m$ be a basis of $H^*(X,{\mathbb{Q}})$, and
let $$\psi_i = c_1(L_i) \in \overline{M}_{g,n}(X,\beta).$$
The {\em descendent} fields, denoted by $\tau_k(\gamma_i)$, correspond 
to the classes $\psi_i^k \text{ev}_i^*(\gamma_j)$ on the moduli space
of maps. 
Let
$$\langle \tau_{k_1}(\gamma_{l_1}) \cdots
\tau_{k_r}(\gamma_{l_r})\rangle_{g,\beta} = \int_{[\overline{M}_{g,r}(X,\beta)]^{vir}} 
\prod_{i=1}^r \psi_i^{k_i} \text{ev}_i^*(\gamma_{l_i})$$
denote the descendent
Gromov-Witten invariants. Foundational aspects of the theory
are treated, for example, in \cite{Beh, BehFan, LiTian}.

Let $C$ be a possibly disconnected curve with at worst nodal singularities.
The genus of $C$ is defined by $1-\chi(\oh_C)$. 
Let $\overline{M}'_{g,r}(X,\beta)$ denote the moduli space of maps
with possibly {disconnected} domain
curves $C$ of genus $g$ with {\em no} collapsed connected components.
The latter condition requires 
 each connected component of $C$ to represent
a nontrivial class in $H_2(X,{\mathbb Z})$. In particular, 
$C$ must represent a {\em nonzero} class $\beta$.

The descendent invariants are defined in the disconnected case by
$$\langle \tau_{k_1}(\gamma_{l_1}) \cdots
\tau_{k_r}(\gamma_{l_r})\rangle'_{g,\beta} = \int_{[\overline{M}'_{g,r}(X,\beta)]^{vir}} 
\prod_{i=1}^r \psi_i^{k_i} \text{ev}_i^*(\gamma_{l_i}).$$
Define the following generating function, 
\begin{equation}
\label{abc}
\bZ'_{GW}\Big(X;u\ |\ \prod_{i=1}^r \tau_{k_i}(\gamma_{l_i})\Big)_\beta = 
\sum_{g\in{\mathbb Z}} \langle \prod_{i=1}^r
\tau_{k_i}(\gamma_{l_i}) \rangle'_{g,\beta} \ u^{2g-2}.
\end{equation}
Since the domain components must map nontrivially, an elementary
argument shows the genus $g$ in the  sum \eqref{abc} is bounded from below. 
The descendent insertions in \eqref{abc} should
match  the  (genus independent) virtual dimension,
$$\text{dim} \ [\overline{M}'_{g,r}(X,\beta)]^{vir} = \int_\beta c_1(T_X) + r.$$
Following the terminology of \cite{MNOP1},
we view \eqref{abc} as a  {\em reduced} partition function.

\subsection{{DT theory}}
\label{uch}
Donaldson-Thomas theory is defined via integration over
the moduli space of ideal sheaves.
Let $X$ be a nonsingular, projective  3-fold.
An {\em ideal sheaf} is a torsion-free sheaf of rank 1 with
trivial determinant. Each ideal sheaf ${\mathcal I}$
injects into its double dual,
$$0 \rarr {\mathcal I} \rarr {\mathcal I} ^{\vee\vee}.$$
As ${\mathcal I}^{\vee\vee}$ is reflexive of rank 1 with trivial
determinant,
$${\mathcal I}^{\vee\vee}\eqq \oh_X,$$
see \cite{OSS}.
Each ideal sheaf $\cI$
determines a subscheme $Y\subset X$,
$$0 \rarr {\mathcal I} \rarr \oh_X \rarr \oh_{Y} \rarr 0.$$
The maximal dimensional components of $Y$
(weighted by their intrinsic multiplicities) determine
an element,
$$[Y] \in H_*(X,\Z).$$
Let $I_n(X,\beta)$ denote the moduli space of ideal sheaves
${\mathcal I}$ satisfying
$$\chi(\oh_{Y}) = n,$$
and
$$[Y] = \beta\in H_2(X,\Z).$$
Here, $\chi$ denotes the holomorphic Euler characteristic.

The Donaldson-Thomas invariant is defined
via integration against virtual class,
$$[I_n(X,\beta)]^{vir}.$$
Foundational aspects of the theory
are treated in \cite{MP,Thom}.

\begin{lm} The virtual dimension of $I_n(X,\beta)$ equals $\int_\beta c_1(T_X)$.
\label{dimm}
\end{lm}
\bpf
The virtual dimension, obtained from the obstruction theory, is
$$\chi(\oh_X,\oh_X) -\chi(\mathcal I, \mathcal I),$$
where 
$$\chi(A,B) = \sum_{i=0}^3  (-1)^{i} \text{dim} \ext^i(A,B).$$
Since $X$ is a nonsingular 3-fold, there exists a finite resolution of
$\mathcal I$ by locally free sheaves,
$$0 \rarr F_3 \rarr F_2 \rarr F_1 \rarr F_0 \rarr \mathcal I \rarr 0.$$
Let $x_{ij}$ denote the Chern roots of $F_i$.
Since the determinant of $\mathcal I$ is trivial, 
$$\sum_{i=0}^3\sum_j (-1)^i x_{ij} = 0.$$
Since the fundamental class of $Y$ is $\beta$,
$$-\ch_2(\mathcal I) = \ch_2(\oh_Y) = \beta.$$

We will calculate the virtual dimension in terms of the
Chern roots via GRR. The first term is,
\begin{equation}
\label{llpp}
\chi(\oh_X,\oh_X)= \int_X \text{Td}(X).
\end{equation}
Next,
$$-\chi(\mathcal I, \mathcal I) = -\int_X \big( \sum_{i=0}^3 \sum_j (-1)^i e^{-x_{ij}}\big) \cdot
\big( \sum_{\hat{i}=0}^3 \sum_{\hat{j}} (-1)^{\hat{i}}e^{x_{\hat{i}\hat{j}}}\big) \cdot \text{Td}(X).$$
Since the Chern root expression in the integrand is {\em even}, only the components
in degrees 0 and 2 need be considered.
The degree 0 component is equal to 1, the square of the rank of $\mathcal I$. The integral
of the degree 0 component against $\text{Td}(X)$
cancels the first term \eqref{llpp}. The degree 2 component is
$$\sum_{i,\hat{i}=0}^3 \sum_{j,\hat{j}}(-1)^{i+\hat{i}} \left(
\frac{x^2_{ij}}{2}-x_{ij}x_{\hat{i}\hat{j}}+\frac{x^2_{\hat{i}\hat{j}}}{2} 
\right) =2\ch_2(\mathcal I) -
\sum_{i,\hat{i}=0}^3 \sum_{j,\hat{j}} (-i)^{i+\hat{i}}x_{ij} x_{\hat{i}\hat{j}}.$$
The second term on the right equals the square of the determinant of $\mathcal I$
and hence vanishes.
We conclude the virtual dimension equals
$$-\int_X 2 \ch_2(\mathcal I) \cdot \text{Td}(X) = \int_\beta c_1(X)$$ 
since the degree 1 term of $\text{Td}(X)$ is $c_1(X)/2$.   
\epf

The moduli space $I_n(X,\beta)$ is canonically isomorphic to the
Hilbert scheme \cite{MP}. As the Hilbert scheme is a fine moduli
space, universal structures are well-defined.
Let $\pi_1$ and $\pi_2$ denote the projections to the respective
factors of $I_n(X,\beta) \times X$.
Consider the universal ideal sheaf $\mathfrak{I}$,
$$\mathfrak{I} \rarr I_n(X,\beta) \times X.$$
Since $\mathfrak{I}$ is $\pi_1$-flat and $X$ is nonsingular, a finite resolution
of $\mathfrak{I}$ by locally free sheaves on $I_n(X,\beta) \times X$
exists. Hence, the Chern classes of $\mathfrak{I}$ are well-defined.

For $\gamma\in H^l(X,\mathbb{Z})$, 
let $\ch_{k+2}(\gamma)$ denote the following operation on the 
homology of $I_n(X,\beta)$:
$$\ch_{k+2}(\gamma):H_*(I_n(X,\beta),\mathbb{Q}) \rarr H_{*-2k+2-l}(I_n(X,\beta),\mathbb{Q}),$$
$$\ch_{k+2}(\gamma)\big(\xi\big)=
 \pi_{1*}\big(   \ch_{k+2}(\mathfrak{I})\cdot \pi_2^*(\gamma) \cap      \pi_1^*(\xi)\big).$$

We define descendent fields in Donaldson-Thomas theory, denoted by
$\tilde{\tau}_k(\gamma)$, to correspond to the operations $(-1)^{k+1}\ch_{k+2}(\gamma)$. 
The descendent invariants are defined by
$$\langle \tilde{\tau}_{k_1}(\gamma_{l_1}) \cdots
\tilde{\tau}_{k_r}(\gamma_{l_r})\rangle_{n,\beta} = \int_{[{I}_n(X,\beta)]^{vir}} 
\prod_{i=1}^r (-1)^{k_i+1} \ch_{k_i+2}(\gamma_{l_i}),$$
where the latter integral is the push-forward to a point of
the class
$$(-1)^{k_1+1}\ch_{k_1+2}(\gamma_{l_1}) \ \circ\ \cdots\
\circ\  (-1)^{k_r+1}\ch_{k_r+2}(\gamma_{l_r})\Big( [I_n(X,\beta)]^{vir}\Big).$$
A similar slant product construction can be found in the Donaldson
theory of 4-manifolds. Since the Chern character contains denominators,
the descendent invariants in Donaldson-Thomas theory
are rational numbers.

Define the Donaldson-Thomas partition function with descendent insertions
by
\begin{equation}
\label{abb}
\bZ_{DT}\Big(X;q\ |\ \prod_{i=1}^r \tilde{\tau}_{k_i}(\gamma_{l_i})\Big)_\beta = 
\sum_{n\in{\mathbb Z}} \langle \prod_{i=1}^r
\tilde{\tau}_{k_i}(\gamma_{l_i}) \rangle_{n,\beta} \ q^{n}.
\end{equation}
An elementary
argument shows the charge $n$ in the  sum \eqref{abb} is bounded from below. 
As before, the descendent insertions in \eqref{abb} should
match  the virtual dimension.

The reduced  partition function is obtained by formally removing the
degree 0 contributions,
$$\bZ'_{DT}\Big(X;q\ |\ \prod_{i=1}^r \tilde{\tau}_{k_i}(\gamma_{l_i})\Big)_\beta =
\frac{\bZ_{DT}\Big(X;q\ |\ \prod_{i=1}^r \tilde{\tau}_{k_i}(\gamma_{l_i})\Big)_\beta}
{\bZ_{DT}(X;q)_0}.$$

The degree 0 partition function is determined by a conjecture of \cite{MNOP1}.
Following the Calabi-Yau case, we conjecture the series $\bZ'_{DT}$ to be a rational
function of $q$.

\vspace{+10pt}
\noindent {{\bf Conjecture 1.}}
The degree 0 Donaldson-Thomas partition function for a 3-fold $X$ is
determined by:
$$
\bZ_{DT}(X;q)_0= M(-q)^{\int_X c_3(T_X \otimes K_X)}\,,
$$
where
$$
M(q) = \prod_{n\geq 1} \frac{1}{(1-q^n)^n}
$$
is the McMahon function.
\vspace{+10pt}

\vspace{+10pt}
\noindent {{\bf Conjecture 2.}}
The reduced series
$\bZ'_{DT}\big(X;q\ | \  \prod_{i=1}^r \tilde{\tau}_{k_i}(\gamma_{l_i})    \big)_\beta$ is
 a rational function of $q$.
\vspace{+10pt}

\subsection{Primary fields}

The GW/DT correspondence is easiest to state for the
{\em primary} fields $\tau_0(\gamma)$ and $\tilde{\tau}_0(\gamma)$.

\vspace{+10pt}
\noindent {{\bf Conjecture 3.}} After the change of variables $e^{iu}=-q$,
$$
(-iu)^d\  \bZ'_{GW}\left(X;u\ | \ \prod_{i=1}^r {\tau}_{0}(\gamma_{l_i})  \right)_{\beta} = 
(-q)^{-d/2}\ \bZ'_{DT}\left(X;q\ | \  \prod_{i=1}^r \tilde{\tau}_{0}(\gamma_{l_i}) 
\right)_{\beta} \, ,
$$
where $d= \int_\beta c_1(T_X)$.
\vspace{+10pt}

Conjecture 3 is consistent with the calculation of degenerate
contributions in \cite{P}.
Let $C$ be a nonsingular, genus $g$ curve in $X$ which rigidly intersects
cycles dual to the classes $\gamma_{l_1}$, \ldots $\gamma_{l_r}$. The local 
Gromov-Witten series is determined in \cite{P},
$$\bZ'_{GW}\left(X;u\ | \ \prod_{i=1}^r {\tau}_{0}(\gamma_{l_i})  \right)_{[C]}=
    \left(\frac{\sin(u/2)}{u/2}\right)^{2g-2 + d}  u^{2g-2},$$ 
The local Donaldson-Thomas series is then {\em predicted} by Conjecture 3,
\begin{eqnarray*}
\bZ'_{DT}\left(X;q\ | \ \prod_{i=1}^r \tilde{\tau}_{0}(\gamma_{l_i})  \right)_{[C]} & = &
   (-iu)^{d}   (-q)^{d/2}   \left(\frac{e^{iu/2}-e^{-iu/2}}{iu}\right)^{2g-2+d} 
             u^{2g-2}  \\
%& = & (-q)^{dim/2} (-q)^{-d/2+1-g}  (-q-1)^{2g-2+d} (i)^{2-2g} (-1)^{d} \\ 
                              & = &     q^{1-g} (1+q)^{2g-2+d}
\end{eqnarray*}
The normalizations and signs in Conjecture 3 are fixed by the requirement that the
reduced partition function
$\bZ'_{DT}$ has initial term $q^{1-g}$ corresponding to the ideal of $C$.

If the cohomology classes $\gamma_i$ are integral, 
the Donaldson-Thomas invariants for primary fields are {\em integer} valued.
The integrality constraints for Gromov-Witten theory obtained via the GW/DT correspondence
for primary fields were conjectured previously in \cite{P,icmtalk}.

\subsection{Descendent fields}
\label{skk}
For fixed curve class $\beta$, consider the full set of (normalized) reduced
partition functions,
$$ {\bZ'}_{GW,\beta}= \left\{ 
(-iu)^{d-\sum
{k_i}}\  \bZ'_{GW}\left(X;u\ | \ \prod {\tau}_{k_i}(\gamma_{l_i})  \right)_{\beta}\right\},$$
where $d=\int_\beta c_1(T_X)$ as before.
Here, ${\bZ'}_{GW,\beta}$ consists of the {\em finite} set of descendent series
with insertions of the correct dimension. 
The 
 set ${\bZ'}_{GW,\beta}$ is 
partially ordered by $\sum k_i$, the {\em descendent} partial ordering.
Similarly, let
$$ {\bZ'}_{DT,\beta}= \left\{ 
(-q)^{-d/2}\  \bZ'_{DT}\left(X;q\ | \ \prod \tilde{\tau}_{k_i}(\gamma_{l_i})
  \right)_{\beta}\right\}.$$

\vspace{+10pt}
\noindent {{\bf Conjecture 4.}} After the change of variables $e^{iu}=-q$, 
\begin{enumerate}
\item[(i)] the sets of functions ${\bZ'}_{GW,\beta}$ and ${\bZ'}_{DT,\beta}$
have the {\em same} linear spans,
\item[(ii)] there exists a canonical matrix expressing the functions ${\bZ'}_{GW,\beta}$
as linear combinations of the functions ${\bZ'}_{DT,\beta}$ :
\begin{enumerate}
\item[(a)] the matrix coefficients depend only upon the classical cohomology of $X$
           and universal series,
\item[(b)] the matrix is unipotent and upper-triangular with respect to the descendent
partial ordering.
\end{enumerate}
\end{enumerate}
\vspace{+10pt}

\noindent By Conjecture 4, each element of ${\bZ'}_{GW,\beta}$ is a canonical linear
combination,
\begin{equation}
\label{jwcv}
(-iu)^{d-\sum
{k_i}}\  \bZ'_{GW}\left( \prod {\tau}_{k_i}(\gamma_{l_i})  \right)_{\beta} =
(-q)^{-d/2}\  \bZ'_{DT}\left( \prod \tilde{\tau}_{k_i}(\gamma_{l_i})  \right)_{\beta} + ...,
\end{equation}
where the omitted terms are strictly lower in the partial ordering.

We do not yet have a complete formula for the canonical matrix of Conjecture 4.
However, for the descendents of the point class $[P]\in H^6(X,\mathbb{Z})$, we
 can formulate a precise
conjecture.

\vspace{+10pt}
\noindent {{\bf Conjecture $4'$.}} 
After the change of variables $e^{iu}=-q$,
\begin{multline*}
(-iu)^{d-\sum k_j}\ 
 \bZ'_{GW}\left( \prod \tau_0(\gamma_{l_i}) \prod {\tau}_{k_j}(P)  \right)_{\beta} 
= \\
(-q)^{-d/2}\ \bZ'_{DT}\left( \prod \tilde{\tau}_0(\gamma_{l_i})
 \prod \tilde{\tau}_{k_j}(P) \right)_{\beta} \, ,
\end{multline*}
if $\text{codim}(\gamma_{l_i})>0$ for each $i$.
\vspace{+10pt}

\subsection{Reactions}
We believe the upper-triangular matrix of Conjecture 4 is determined by
two types of {\em reactions}:
\begin{eqnarray*}
\tau_a(\gamma_l) & \rarr &  A_{a}(\gamma_l)  \ \tau_{a-1}(c_1(T_X) \cup \gamma_l) \\
\tau_a(\gamma_l) {\tau}_{a'}(\gamma_{l'})  & \rarr &  
A_{a,a'}(\gamma_l,\gamma_{l'}) \ \tau_{a+a'-1}(\gamma_l \cup
\gamma_l') 
\end{eqnarray*}
The linear combination \eqref{jwcv} should be generated by applying the two 
reactions to the Gromov-Witten insertions $$\prod \tau_{k_i}(\gamma_{l_i})$$
to exhaustion and then interpreting  the output in Donaldson-Thomas theory. For example,
\begin{multline*}
(-iu)^{d-
{k}}\  \bZ'_{GW}\left( {\tau}_{k}(\gamma_{l})\right)_{\beta} = \\
(-q)^{-d/2}
\sum_{j=0}^k \left(\prod_{i=1}^j A_{k-i+1}(c_1(T_X)^{i-1} \cup \gamma_l)\right) \ 
\  \bZ'_{DT}\left( \tilde{\tau}_{k-j}(c_1(T_X)^j \cup \gamma_{l})\right).
\end{multline*}
The resulting matrix will be upper-triangular with respect to the
{\em reaction partial ordering}, a refinement of the descendent
partial ordering.
We further speculate that the {\em reaction amplitudes},
$$A_{a}(\gamma_l),\  A_{a,a'}(\gamma_l,\gamma_{l'}),$$
are given by universal formulas depending only upon  the classical cohomology of
$X$ (including possibly the Hodge decomposition).
Conjectures 3, 4, and $4'$ are all {\em consequences} of the reaction view of the
GW/DT correspondence for descendent fields.

\subsection{An example}
Let $X$ be $\proj^3$ and let $\beta$ be the class $[L]$ of a line. A Gromov-Witten calculation
using localization and known Hodge integral evaluations yields the following
result,
$$\bZ'_{GW}\Big(X;u\ | \ {\tau}_{0}(L) \tau_1(P)  \Big)_{[L]}= 
       \left(\frac{\sin(u/2)}{u/2}\right) \cos(u/2)  u^{-2},$$
see \cite{GrabPan,FabPan}.
By Conjecture $4'$,
\begin{eqnarray*}
\bZ'_{DT} \Big(X;q\ |\ \tilde{\tau}_{0}(L) \tilde{\tau}_{1}(P)  \Big)_{[L]} & = & 
(-iu)^3 (-q)^2 \Big(\frac{\sin(u/2)}{u/2}\Big) \cos(u/2) u^{-2} \\
& = & (-iu)^3 (-q)^2 \frac{ e^{iu/2}-e^{-iu/2}}{iu} \frac{e^{iu/2}+e^{-iu/2}}{2} u^{-2}\\
& = & \frac{1}{2} q (1-q^2)
\end{eqnarray*}
The resulting Donaldson-Thomas series can be checked order by order in $q$
via localization.

\section{The GW/DT correspondence for relative theories}
\subsection{GW theory}
Let $X$ be a nonsingular, projective 3-fold and let $S\subset X$ be a
nonsingular divisor.
The Gromov-Witten theory of $X$ relative to $S$ has been defined in
\cite{EGH,IP,LR,L}.
Let $\beta\in H_2(X,{\mathbb Z})$ be a curve class satisfying
$$\int_\beta [S] \geq 0.$$
Let $\stackrel{\rightarrow}{\mu}$ be an ordered partition,
$$\sum \mu_j = \int_\beta[S],$$
with positive parts.
The moduli space $\overline{M}'_{g,n}(X/S,\beta, \stackrel{\rightarrow}{\mu})$ parameterizes
stable relative maps with possibly disconnected domains and
relative multiplicities determined by $\stackrel{\rightarrow}{\mu}$.
 As usual, 
the connected components of the domain  are required to
map nontrivially. The target of a relative map is allowed to be
a $k$-step degeneration, $X[k]$, of $X$ along $S$, see \cite{L}.

The relative conditions in the theory correspond to partitions {\em weighted} by
the cohomology of $S$.
Let  $\delta_1, \ldots, \delta_{m_S}$ be a basis of $H^*(S,{\mathbb Q})$.
A cohomology weighted partition $\mathbf{\eta}$
consists of an {\em unordered} set of pairs,
$$\left\{ (\eta_1, \delta_{\ell_1}), \ldots, (\eta_{s}, \delta_{\ell_{s}}) 
\right\}, $$
where $\sum_j \eta_j$ is an {\em unordered} partition of $\int_\beta [S]$.
The automorphism group, $\text{Aut}(\mathbf{\eta})$, consists of 
permutation symmetries of ${\mathbf \eta}$.

The {\em standard} order on the parts of $\eta$ is
$$(\eta_i,\delta_{\ell_i})> (\eta_{i'},\delta_{\ell_{i'}})$$
if $\eta_i>\eta_{i'}$ or if $\eta_i=\eta_{i'}$ and $\ell_i>\ell_{i'}$.
Let $\stackrel{\rarr}{\eta}$ denote the partition $(\eta_1,\ldots, \eta_{s})$
obtained from the standard order.

Relative Gromov-Witten invariants are defined by integration against the
virtual class of the moduli of maps. 
Let $\gamma_1,\ldots,\gamma_{m_X}$ be a basis of $H^*(X,{\mathbb Q})$,
and let
\begin{multline*}
\langle \tau_{k_1}(\gamma_{l_1}) \cdots
\tau_{k_r}(\gamma_{l_r})\ | {\mathbf \eta} 
\rangle'_{g,\beta} = \\
 \frac{1}{|\text{Aut}(\mathbf \eta)|}
\int_{[\overline{M}'_{g,r}(X/S,\beta,\stackrel{\rarr}{\eta})]^{vir}} 
\prod_{i=1}^r \psi_i^{k_i}\text{ev}_i^*(\gamma_{l_i}) \cup 
\prod_{j=1}^{s} \text{ev}^*_j(\delta_{\ell_j}).
\end{multline*}
Here, the second evaluations, 
$$\text{ev}_j: \overline{M}'_{g,r}(X/S,\beta,\stackrel{\rarr}\eta) \rarr S.$$
are determined by the relative points.

The  Gromov-Witten invariant is defined for {\em unordered}
weighted partitions $\mathbf{\eta}$. However, to fix the sign, the integrand
on the right side requires an ordering. The ordering is corrected
by the automorphism prefactor.

As before, we will require the associated Gromov-Witten partition function, 
\begin{equation}
\label{ab}
\bZ'_{GW}\Big(X/S ; u\ |\ \prod_{i=1}^r \tau_{k_i}(\gamma_{l_i}) \Big)_{\beta,{\mathbf \eta}}
 = 
\sum_{g\in{\mathbb Z}} \langle \prod_{i=1}^r
\tau_{k_i}(\gamma_{l_i}) \ | {\mathbf \eta} \rangle'_{g,\beta} \ u^{2g-2}.
\end{equation}
The definitions here parallel those of Section \ref{gwt}.

\subsection{DT theory}
\subsubsection{Stable relative ideal sheaves}
\label{relid}
Relative Donaldson-Thomas theory is defined via integration over the moduli
space of relative ideal sheaves. We outline J. Li's definition
of the relative theory here \cite{jli}. A full foundational treatment has not yet been
written.

Let $X$ be a nonsingular, projective 3-fold and let $S\subset X$ be a
nonsingular divisor.
Let $\cI$ be an ideal sheaf on $X$ with associated subscheme $Y$.
The ideal sheaf $\cI$ is {\em relative} to $S$ if the natural map,
$$\cI\otimes_{\oh_X}\oh_S \rarr \oh_X \otimes_{\oh_X}\cO_S,
$$
is injective. Relativity may be viewed as a transversality condition of $Y$ with
respect to $S$.
In particular, the scheme theoretic intersection,
$Y \cap S,$
defines an element of  the Hilbert scheme, $$\text{Hilb}(S, \int_\beta [S]),$$
 of points of $S$.

Relativity is an {\em open} condition on ideal sheaves on $X$. 
A proper moduli
space, $I_n(X/S,\beta)$, of relative ideal sheaves is constructed by
considering stable ideal sheaves relative on 
the degenerations $X[k]$ of $X$.

Let $S_0, \ldots, S_k$ denote the canonical images of $S$ in the 
degeneration $X[k]$. Here, $S_0, \ldots, S_{k-1}$ are the singular
divisors, and $S_k$ is the transform of the original relative divisor.
An ideal sheaf on $X[k]$ is {\em predeformable} if,
for every singular divisor $S_l\subset X[k]$, the induced map,
$$\cI\otimes_{\oh_{X[k]}}\oh_{S_l}\rarr \oh_{X[k]}\otimes_{\oh_{X[n]}}\oh_{S_l}$$
is injective. 

Let $Y_0,\ldots, Y_k$ be the restrictions of $Y$ to the
components of $X[k]$ with $Y_{k}$ and $Y_{k+1}$ incident $S_k$. 
The predeformability condition at the singular divisor $S_l$ can be
restated in the following form:
$Y_l$ and $Y_{l+1}$ are transverse to $S_l$ with
equal
scheme theoretic intersections,
$$Y_l \cap S_l = Y_{l+1} \cap S_l \subset S_l.$$

Ideal sheaves $\cI_1$ and
$\cI_2$ on the degenerations $X[k_1]$ and $X[k_2]$ are isomorphic
 if $k_1=k_2$ and there exists an isomorphism of varieties $$\sigma:
X[k_1]\rarr X[k_2]$$ over $X$ such that
$$
\sigma^* \{\cI_2\to\oh_{X[k_2]}\}\cong \{\cI_1\to\oh_{X[k_1]}\},
$$
where the isomorphism $\sigma^* \oh_{X[k_2]}\cong \oh_{X[k_1]}$
is the identity. The automorphism group, $\text{Aut}(\cI)$, is
the set of equivalences of $\cI$ to itself. 
A predeformable ideal sheaf $\cI$ on $X[k]$ relative to $S_k$ is
{\em stable} if $\text{Aut}(\cI)$ is finite.

The moduli space, $I_n(X/S,\beta)$, parameterizes stable, predeformable,
ideal sheaves $\cI$ on degenerations $X[k]$ relative $S_k$ satisfying 
$$\chi(\oh_Y)=n$$ and
$$\pi_*[Y] = \beta \in H_2(X,\Z),$$
where $\pi: X[k] \rarr X$ is the canonical stabilization map.
The moduli space $I_n(X/S,\beta)$  is a complete, Deligne-Mumford stack
equipped with a canonical perfect obstruction theory. 

Relative
Donaldson-Thomas theory is defined via integration against
the associated virtual class. The primary and descendent fields
are  
defined via the Chern characters of the universal ideal
sheaf $\mathfrak{I}$ on the universal product 
stack following Section \ref{uch}. The predeformability
condition is expected to imply the existence of finite resolutions of
$\mathfrak{J}$ by locally free sheaves. The relative
conditions in the theory are defined via the canonical
intersection map,
$$\epsilon: I_n(X/S,\beta) \rarr \text{Hilb}(S,\int_\beta[S]),$$
to the Hilbert scheme of points.

\subsubsection{The Nakajima basis}
The cohomology of the Hilbert scheme of points
of $S$ has a canonical basis indexed by
cohomology weighted partitions.
The basis is obtained from the representation
of the Heisenberg algebra on the cohomologies of the Hilbert schemes of points \cite{Gron, Nak}.

Let ${\mathbf \eta}$ be a cohomology weighted partition 
with respect to the basis 
$\delta_1, \ldots, \delta_{m_S}$
of $H^*(S,\Q)$. 
Following the notation of \cite{Nak},
let
\begin{equation}\label{hdeff}
C_{\mathbf\eta}= \frac{1}{\mathfrak{z}(\mathbf \eta)}\  P_{\delta_1}[\eta_1]
\cdots P_{\delta_s}[\eta_s]\cdot {\mathbf 1}\ \in H^*(\text{Hilb}(S,|\eta|), \Q),
\end{equation}
where
$$
\mathfrak{z}({\mathbf \eta})= {\prod_i \eta_i \ |\text{Aut}({\mathbf \eta})|}\ ,$$
and $|{\mathbf \eta}|= \sum_j \eta_j$.
In the presence of odd cohomology, the sign of
$C_\eta$ is fixed by placing the operator product \eqref{hdeff} in
standard order.

The {\em Nakajima basis} of the cohomology of $\text{Hilb}(S,k)$ is
the set,
$$\left\{ C_{\mathbf \eta} \right \}_{|{\mathbf \eta}|=k},$$
see \cite{Nak}.

We assume
the cohomology basis of $S$ is self dual with respect to the
Poincar\'e pairing. 
Then, to each weighted partition ${\mathbf \eta}$, a dual partition
${\mathbf \eta}^\vee$ is defined by taking the Poincar\'e duals
of the cohomology weights. The Nakajima basis is orthogonal
with respect to the Poincar\'e pairing on the 
cohomology of the Hilbert scheme,
\begin{equation}
\label{ppp}
\int_{\text{Hilb}(S,k)}  C_\eta \cup C_{\nu} =
\frac{(-1)^{k-\ell({\mathbf \eta})}}{\mathfrak{z}(\mathbf{\eta})} \ 
\delta_{ {\mathbf {\nu}}, {\mathbf {\eta}}^\vee},
\end{equation}
see \cite{ES,Nak}.

\subsubsection{Relative Donaldson-Thomas invariants}
Relative Donaldson-Thomas invariants are defined via integration over
the moduli spaces of stable relative sheaves.
The relative and absolute moduli spaces have the same 
deformation theories on the open set of ideal sheaves on $X$ relative to $S$. Hence,
$I_n(X/S,\beta)$ has virtual dimension $\int_\beta c_1(X)$ by Lemma \ref{dimm}.
Though the open set may be {\em empty}, the conclusion is correct, and
we leave the details to the reader.

The descendent invariants in relative Donaldson-Thomas theory are defined by
$$\langle \tilde{\tau}_{k_1}(\gamma_{l_1}) \cdots
\tilde{\tau}_{k_r}(\gamma_{l_r})\ | {\mathbf{\eta}}
\rangle_{n,\beta} = \int_{[{I}_n(X/S,\beta)]^{vir}} 
\prod_{i=1}^r (-1)^{k_i+1} \ch_{k_i+2}(\gamma_{l_i}) \cap \epsilon^*(C_{\mathbf \eta}).$$
Define the associated partition function by
\begin{equation}
\label{abbb}
\bZ_{DT}\Big(X/S;q\ |\ \prod_{i=1}^r \tilde{\tau}_{k_i}(\gamma_{l_i})\Big)_{\beta, {\mathbf \eta}} 
= 
\sum_{n\in{\mathbb Z}} \langle \prod_{i=1}^r
\tilde{\tau}_{k_i}(\gamma_{l_i})\ | {\mathbf \eta} \rangle_{n,\beta} \ q^{n}.
\end{equation}
As before the charge $n$ in the  sum \eqref{abb} is bounded from below. 

The reduced  partition function is obtained by formally removing the
degree 0 contributions,
$$\bZ'_{DT}\Big(X/S;q\ |\ \prod_{i=1}^r \tilde{\tau}_{k_i}(\gamma_{l_i})\Big)_
{\beta, {\mathbf \eta}} =
\frac{\bZ_{DT}\Big(X/S;q\ |\ \prod_{i=1}^r \tilde{\tau}_{k_i}(\gamma_{l_i})
\Big)_{\beta, {\mathbf \eta}}}
{\bZ_{DT}(X/S;q)_0}.$$

We conjecture a complete formula for degree 0 relative theory.
Let $\Omega_X[S]$ denote the locally free sheaf of
differential forms of $X$ with logarithmic poles
along $S$. Let 
$$T_X[-S]= \Omega_X[S]\ ^\vee,$$
denote the dual sheaf of tangent fields with logarithmic zeros.
Let 
$$K_X[S]= \Lambda^3 \Omega_X[S]$$
denote the logarithmic canonical class.

\vspace{+10pt}
\noindent {{\bf Conjecture 1R.}}
The degree 0 relative Donaldson-Thomas partition function for a 3-fold $X$ is
determined by:
$$
\bZ_{DT}(X/S;q)_0= M(-q)^{\int_X c_3(T_X[-S] \otimes K_X[S])}\,.
$$
\vspace{+10pt}

If $S$ is empty, Conjecture 1R specializes to Conjecture $1'$ of \cite{MNOP1}.
A proof of Conjecture 1R in the toric case is presented in
Section \ref{evv}.
As before, we conjecture the reduced series are
rational functions of $q$.

\vspace{+10pt}
\noindent {{\bf Conjecture 2R.}}
The reduced series
$\bZ'_{DT}\big(X/S;q\ | \  \prod_{i=1}^r \tilde{\tau}_{k_i}(\gamma_{l_i})   \big)_{
\beta, {\mathbf \eta}}$ is
 a rational function of $q$.
\vspace{+10pt}

\subsection{Primary fields}

We restrict our discussion of the relative GW/DT correspondence to the
primary fields. A treatment of the descendent correspondence
at the level of Section \ref{skk} is left to the reader. In particular, we
do not know the precise formulas for the descendent correspondence.

\vspace{+10pt}
\noindent {{\bf Conjecture 3R.}} After the change of variables $e^{iu}=-q$,
\begin{multline*}
(-iu)^{d+ \ell({\mathbf \eta})-|{\mathbf \eta}|}  
\bZ'_{GW}\left(X/S;u\ | \ \prod_{i=1}^r {\tau}_{0}(\gamma_{l_i})  \right)_
{\beta,{\mathbf \eta}} =  \\
(-q)^{-d/2}\ \bZ'_{DT}\left(X/S;q\ | \  \prod_{i=1}^r \tilde{\tau}_{0}(\gamma_{l_i}) \right)_
{\beta, {\mathbf \eta}} \, ,
\end{multline*}
where $d= \int_\beta c_1(T_X)$ and $\ell({\mathbf \eta})$ denotes the length of
${\mathbf \eta}$.
\vspace{+10pt}

We present the
simplest example in which all the features of the correspondence are visible.
Let $D$ be a nonsingular surface, and let
 $$X=\proj^1 \times D.$$ Let $0, \infty \in  \proj^1$ be two points in the base, and
let $D_0$ and $D_\infty$ be the associated fibers. Let $$[\proj^1]\in H_2(X,\Z)$$
denote the class of the horizontal $\proj^1$, and let 
$\beta=m[\proj^1]$.
We will consider the
theories of $X$ relative to the divisors $D_0$ and $D_\infty$ in the curve class
$\beta$.

Since there are two divisors, the boundary conditions
of the relative theories are specified by two 
partitions $\mathbf \eta$ and $\mathbf \nu$ weighted by the
cohomology of $D$.
The relative Gromov-Witten theory is particularly simple to compute. A direct
calculation
yields the answer,
$$
\bZ'_{GW}\left(X/S;u \right)_
{\beta,{\mathbf \eta}, {\mathbf \nu}} =  \frac{1}{\mathfrak{z}(\mathbf{\eta})}
u^{-2\ell({\mathbf {\eta}})} \ \delta_{ {\mathbf {\nu}}, {\mathbf {\eta}}^\vee}.$$
Our correspondence {\em predicts} the associated Donaldson-Thomas series,
\begin{eqnarray*}
\bZ'_{DT}\left(X/S;q \right)_
{\beta,{\mathbf \eta}, {\mathbf \nu}}& =& (-q)^{d/2} (-iu)^{d-2m+\ell({\mathbf \eta})
+\ell({\mathbf \nu})} \frac{1}{\mathfrak{z}(\mathbf{\eta})}
u^{-2\ell({\mathbf {\eta}})} \ \delta_{ {\mathbf {\nu}}, {\mathbf {\eta}}^\vee} \\
 & = & \frac{(-1)^{m-\ell({\mathbf \eta})}}{\mathfrak{z}(\mathbf{\eta})} q^m \ 
\delta_{ {\mathbf {\nu}}, {\mathbf {\eta}}^\vee} ,
\end{eqnarray*}
using the relation $d=2m$ in the last equality.

The moduli space $I_m(X/D_0\cup D_\infty, \beta)$ is isomorphic to 
$\text{Hilb}(D,m)$. The Donaldson-Thomas invariant is therefore
a classical intersection product,
$$\langle {\mathbf \eta}|\ | {\mathcal \nu} \rangle_{m,\beta} =
\int_{\text{Hilb}(D,m)} C_{\mathbf \eta} \cup C_{\mathbf \nu}.$$
The $q^m$ term of the predicted Donaldson-Thomas series is thus
correct by \eqref{ppp}. The division of the degree 0 series
does {\em not} affect the first term.

\subsection{The degeneration formula}
The relative theories satisfy degeneration formulas. Let 
$$\lambda:{\mathcal X} \rarr C$$
be a nonsingular 4-fold fibered over a nonsingular, irreducible curve.
Let $X$ be a nonsingular fiber of $\lambda$, and
let $$X_1 \cup_S X_2$$
be a reducible special fiber consisting of two 
nonsingular 3-folds intersecting transversely along a nonsingular
surface $S$. The degeneration formulas express the absolute invariants
of $X$ via the relative invariants of $X_1/S$ and $X_2/S$. We will
show the degeneration formulas of the relative theories are
compatible with the GW/DT correspondence for primary fields.

The degeneration formula for Gromov-Witten theory is naturally written
in terms of the absolute and relative partition functions,
\begin{multline*}
\bZ'_{GW}\left(X|  \prod_{i=1}^r {\tau}_{0}(\gamma_{l_i})  \right)_\beta= \\
%\sum_{\beta_1+\beta_2=\beta} \sum_{I_1\cup I_2=[r]}
\sum \ \bZ'_{GW}\left(\frac{X_1}{S}| \prod_{i\in P_1} {\tau}_{0}(\gamma_{l_i})  \right)_
{\beta_1,{\mathbf \eta}} 
{{\mathfrak z}({\mathbf \eta})}{ u^{2\ell({\mathbf \eta})}} \
\bZ'_{GW}\left(\frac{X_2}{S}| \prod_{i\in P_2} {\tau}_{0}(\gamma_{l_i})  \right)_
{\beta_2,{\mathbf \eta}^\vee},
\end{multline*}
where the sum is over curve splittings $\beta_1+\beta_2=\beta$, 
marking partitions $$P_1 \cup P_2=\{1,\ldots, r\},$$ and 
cohomology weighted partitions ${\mathbf \eta}$.
The central factor on the right accounts for the multiplicities and
the shift in the genus variable $u$. A proof can be found in
\cite{EGH,IP,LR,L}.

The degeneration formula for Donaldson-Thomas theory takes a very similar
form,
\begin{multline*}
\bZ'_{DT}\left(X|  \prod_{i=1}^r \tilde{\tau}_{0}(\gamma_{l_i})  \right)_\beta= \\
%\sum_{\beta_1+\beta_2=\beta} \sum_{I_1\cup I_2=[r]}
\sum  \bZ'_{DT}\left(\frac{X_1}{S}| \prod_{i\in P_1} \tilde{\tau}_{0}(\gamma_{l_i})  \right)_
{\beta_1,{\mathbf \eta}} 
\frac{(-1)^{|\eta|-\ell(\mathbf \eta)}{\mathfrak z}({\mathbf \eta})}{ q^{|{\eta}|}} 
\ \bZ'_{DT}\left(\frac{X_2}{S}| {\prod_{i\in P_2}} \tilde{\tau}_{0}(\gamma_{l_i})  \right)_
{\beta_2,{\mathbf \eta}^\vee},
\end{multline*}
where the sum is as before.
The central factor on the right accounts for the diagonal splitting,
$$[\bigtriangleup] = \sum_{|\eta|=k} (-1)^{k-\ell(\mathbf \eta)}  {{\mathfrak z}({\mathbf \eta})}\
 C_{\mathbf \eta} \otimes C_{{\mathbf \eta}^\vee}
 \ \ \in H^*(\text{Hilb}(S,k) \times
\text{Hilb}(S,k), \Q),$$
 and
the shift in the charge variable $q$. The proof should follow \cite{L} but has
yet to be written.

The compatibility between the degeneration formulas and the GW/DT correspondence
is straightforward. Let $d=\int_\beta c_1(X)$ as before, and let
$$d_i=\int_{\beta_1} c_1(X_i).$$ We have a partition of the total
degree $d$, 
\begin{eqnarray*}
d & = & d_1 +d_2- 2\int_{\beta_1}[S]  \\
   & = & (d_1 - |{\mathbf \eta}| ) + (d_2 - |{\mathbf \eta}^\vee|).
\end{eqnarray*}
Using the degree partition, the degeneration formulas for
the relative theories are equivalent via the GW/DT correspondence.

\section{The equivariant vertex measure}
\label{evv}
\subsection{Summary}

Let $\T$ be a 3-dimensional complex torus
with coordinates $t_i$. Let $\T$ act on $\A^3$ with coordinates
$x_i$ by 
\begin{equation}
  \label{st_act}
  (t_1,t_2,t_3)\cdot x_i = t_i x_i  \,.
\end{equation}
In these coordinates,
the tangent representation at the origin $0\in \A^3$
has character $t^{-1}_1 +t_2^{-1}+ t_3^{-1}$.

Let
$\pi$ be a 3-dimensional partition with three outgoing
2-dimensional partitions $\lambda_1$, $\lambda_2$, and $\lambda_3.$
The equivariant vertex $\bV_\pi$ arises in the
localization formula for the Donaldson-Thomas theory of toric
3-folds \cite{MNOP1}. 

The equivariant vertex determines a natural
$3$-parametric family of measures $\bw$ on 3-dimensional
partitions. The measure of $\pi$ is defined by 
$$
\bw(\pi) = \prod_{k\in \Z^3} \left( s,k \right)^{-\bv_k} \,,
$$
where $s=(s_1,s_2,s_3)$ are parameters, $(\, \cdot\, ,\, \cdot\,)$
denotes the standard inner product, and $\bv_k$ is the coefficient
of $t^k$ in $\bV_\pi$.

Consider the generating series of the equivariant vertex measures of
3-dimensional partitions $\pi$ with fixed outgoing 2-dimensional partitions, 
$$\bW(\lambda_1,\lambda_2,\lambda_3) =
\sum_\pi \bw(\pi) q^{|\pi|}.$$
Here $|\pi|$ is defined as the (signed) number of boxes obtained
by formally removing the infinite outgoing cylinders \cite{MNOP1}.

\begin{tm}
\label{dz}
For finite 3-dimensional partitions,
$$\bW(\emptyset, \emptyset, \emptyset)= M(-q)^{- \frac{(s_1+s_2)(s_1+s_3)(s_2+s_3)}{s_1s_2s_3}}.$$
\end{tm}

%\begin{tm}
%\label{d1}
%For 3-dimensional partitions with a single outgoing column,
%$$\bW(1, \emptyset, \emptyset)= M(-q)^{- \frac{(s_1+s_2)(s_1+s_3)(s_2+s_3)}{s_1s_2s_3}}
% (1-q)^{\frac{s_2+s_3}{s_1}}.$$
%\end{tm}

Our proof is  independent of the conjectural GW/DT correspondence. However, 
relative Donaldson-Thomas theory plays an essential role.

\subsection{Equivariant Donaldson-Thomas theory} 
Let the 1-dimensional torus $\T^1$ act on $\proj^1$ with tangent weights
$-s_1$ and $s_1$ at the fixed points $0$ and $\infty$.
Let the 2-dimensional torus $\T^2$ act on $\com^2$ with weights $-s_2$ and $-s_3$.
The torus $\T=\T^1 \times \T^2$ acts on
$$X=\proj^1 \times \com^2 $$
preserving the divisor $S$ over $\infty$.
We will study the equivariant Donaldson-Thomas theory of $X$ relative to $S$.

Since $X$ is not projective, the non-equivariant
theory
is not well-defined. However, the $\T$-equivariant theory can be defined via the 
 residue since the $\T$-fixed locus of $I_n(X/S,0)$ is projective.
Let $\bZ^{\T}_{DT}(X/S;q)_0$ denote the degree 0 partition function
for the equivariant relative theory.

A rational function $f\in \Q(s_1,s_2,s_3)$ has only {\em monomial poles} in
the variables $s_2$ and $s_3$
if 
$$f(s_1,s_2,s_3) = \frac{p(s_1,s_2,s_3)}{s_2^{k_2} s_3^{k_3}}$$
for $p\in \Q[s_1,s_2,s_3]$ and $k_2,k_3\in \Z$.

\begin{lm}
\label{gwx}
The $q$ coefficients of 
$\bZ^{\T}_{DT}(X/S;q)_0$ have only monomial poles in the variables $s_2$ and $s_3$.
\end{lm}

\bpf
The Hilbert-Chow morphism and the collapsing
maps,
$$X[k] \rarr X,$$
together yield a $\T$-equivariant, proper morphism,
$$\iota_1:I_n(X/S,0) \rarr \text{Sym}^n(X).$$
The projection $X \rarr \com^2$ yields a
$\T$-equivariant, proper morphism,
$$\iota_2: \text{Sym}^n(X) \rarr \text{Sym}^n(\com^2).$$
Finally, a $\T$-equivariant, proper
morphism,
$$\iota_3: \text{Sym}^n(\com^2) \rarr \oplus_{1}^n\com^2,$$
is obtain via the higher moments,
$$\iota_3\Big( \ \{(x_i,y_i)\} \ \Big) = \Big(\sum_i x_i, \sum_i y_i\Big) \oplus 
\Big(\sum_i x^2_i, \sum_i y^2_i\Big) \oplus \cdots \oplus 
\Big(\sum_i x^n_i, \sum_i y^n_i\Big).$$
Let $j=\iota_3\circ \iota_2\circ \iota_1$.

The virtual class $[I_n(X/S,0)]^{vir}$ is an element of the $\T$-equivariant
Chow ring of $I_n(X/S,0)$. Since $j$ is $\T$-equivariant and proper,
$$\int_{[I_n(X/S,0)]^{vir}} 1 =
\int_{\oplus_1^n \com^2} j_*[I_n(X/S,0)]^{vir},$$
where $\int$ denotes $\T$-equivariant integration.
Since the space $\oplus_1^n \com^2$ has a unique $\T$-fixed point
with tangent weights,
$$-s_2,-s_3,-2s_2,-2s_3, \ldots, -ns_2, -ns_3,$$
we conclude the integral
$$\int_{\oplus_1^n \com^2} j_*[I_n(X/S,0)]^{vir},$$
has only monomial poles in the variables $s_2$ and $s_3$.
\epf

\subsection{Localization}
The $\T$-fixed loci of $I_n(X/S,0)$ lie over either 0 or $\infty$.
The fixed points over $0$ correspond to finite 3-dimensional partitions
with localization contributions to $\bZ^{\T}_{DT}(X/S;q)_0$
determined by $\bW(\emptyset,\emptyset,\emptyset)$, see \cite{MNOP1}.

A Donaldson-Thomas theory of {\em rubber} naturally arises on the 
fixed loci of $I_n(X/S,0)$ over $\infty$.
Let 
$$R= \proj^1 \times \com^2,$$
and let $S_0$ and $S_\infty$ denote the divisors over $0$ and $\infty$
respectively. Let $\T^2$ act on $\com^2$ with weights $-s_2$ and $-s_3$.
We will consider the $\T^2$-equivariant Donaldson-Thomas rubber theory
of $R$ relative to $S_0$ and $S_\infty$. For the rubber theory,
sheaves differing by the $\T^1$ action on $\proj^1$ are identified.
We denote the rubber theory by a superscripted tilde.

The rubber moduli space $I_n(R/S_0\cup S_\infty, 0)\ \tilde{}$ carries
cotangent lines at the dynamical points $0$ and $\infty$ of $\proj^1$. Let
$\psi_0$ denote the class of the cotangent line at $0$.
Let
$$\bW_\infty = 1+\sum_{n\geq 1} \int_{[I_n(R/S_0\cup S_\infty, 0)\ \tilde{}\ ]^{vir}}
\frac{1}{s_1 -\psi_0},$$
where $\int$ here denotes $\T^2$-equivariant integration.
The leading term 1 may be viewed as a degenerate $n=0$ contribution.
By the virtual localization formula applied to the relative Donaldson-Thomas
theory of $X/S$, the series $\bW_\infty$ generates the localization contributions 
to $\bZ^{\T}_{DT}(X/S;q)_0$ of
the $\T$-fixed points over $\infty$.

The product of the localization contributions over $0$ and $\infty$
yields the partition function,
\begin{equation}
\label{pwq}
\bZ^{\T}_{DT}(X/S;q)_0 = \bW(\emptyset,\emptyset,\emptyset)  \cdot \bW_\infty.
\end{equation}

Consider the $\T^2$-equivariant rubber theory {\em without} any cotangent
line insertions,
$$\bF_\infty = \sum_{n\geq 0} q^n \int_{[I_n(R/S_0\cup S_\infty, 0)\ \tilde{}\ ]^{vir}} 1.$$
By definition,
$$\bW(\emptyset,\emptyset,\emptyset), \ \bW_\infty \in \Q(s_1,s_2,s_3)[[q]],$$
 and 
$$\bF_\infty \in \Q(s_2,s_3)[[q]].$$

\begin{lm}
\label{llggf}
$\log W_\infty = \frac{1}{s_1} \bF_\infty$. 
\end{lm}

\bpf
We first expand $\bW_\infty$ by powers of the cotangent line,
$$\bW_\infty = 1+\sum_{l\geq 0} \frac{1}{s_1^{l+1}} \bF_{\infty,l}\, ,$$
where
$$\bF_{\infty,l} = \sum_{n\geq 1} q^n 
\int_{[I_n(R/S_0\cup S_\infty, 0)\ \tilde{}\ ]^{vir}} \psi_0^l.$$

Next, we apply a version of the topological recursion relation
to inductively calculate $\bF_{\infty,l}$.
Let 
$$\pi:{\mathcal Y}_n \rarr I_n(R/S_0\cup S_\infty, 0)\ \tilde{}$$
be the universal subscheme over the moduli space.
The morphism $\pi$ is finite, flat, and compatible with the
$\T^2$-action. Therefore,
$$q\frac{d}{dq} \bF_{\infty,l} =
\sum_{n> 0} q^n 
\int_{[{\mathcal Y}_n]^{vir}} \psi_0^l,$$
where the virtual class of ${\mathcal Y}$ is defined as the pull-back of the
virtual class of the moduli space by $\pi$.
The canonical map,
$$f: {\mathcal Y}_n \rarr R[k],$$
projects further to $\proj^1[k]$, the associated degeneration
of the base $\proj^1$. By the definition of the relative moduli
space, the image in $\proj^1[k]$ is always disjoint from
the relative points $0$ and $\infty$ and the nodes.
Hence, the family of degenerating bases over ${\mathcal Y}_n$
has {\em three} disjoint nonsingular sections.

The application of the topological recursion relation determined
by the three sections to $\psi_0$ yields the following equation,
$$q\frac{d}{dq} \bF_{\infty,l}= \bF_{\infty,l-1} \cdot
q \frac{d}{dq}\bF_{\infty,0}.$$ The solution,
$$\bF_{\infty,l} = \frac{ \bF_{\infty,0}^{l+1}}{(l+1)!},$$
is easily found.
We conclude $\bW_\infty= \exp(\frac{1}{s_1} \bF_\infty)$.
\epf

\subsection{Proof of Theorem \ref{dz}}
The logarithm of equation \eqref{pwq} yields the relation,
 $$\log \bW(\emptyset,\emptyset,\emptyset) = \log \bZ^{\T}_{DT}(X/S;q)
- \log \bW_\infty.$$
By Lemmas \ref{gwx} and \ref{llggf}, the $q$ coefficients of $\log \bW(\emptyset,
\emptyset,\emptyset)$ are of the form
$$\frac{1}{s_1} \frac{p_1(s_1,s_2,s_3)}{p_2(s_2,s_3)},$$
where the $p_i$ are polynomials. Since the equivariant vertex measure is
a degree 0 rational function \cite{MNOP1}, 
$$\text{deg}(p_1) = 1 + \text{deg}(p_2).$$
Since the series $\bW(\emptyset,\emptyset,\emptyset)$ is {\em symmetric} in the
variables $s_1$, $s_2$, and $s_3$, we conclude, 
$$\log \bW(\emptyset,\emptyset,\emptyset) = \frac{1}{s_1s_2s_3} \bF_0(q,s_1,s_2,s_3),$$
where $\bF_0 \in \Q[s_1,s_2,s_3][[q]].$
The coefficients of $\bF_0$ must be {\em cubic} polynomials.

By Lemma \ref{divv} below, the $q^{n}$ coefficient of $\bW(\emptyset,\emptyset,\emptyset)$
is divisible by the cubic factor $(s_1+s_2)(s_1+s_3)(s_2+s_3)$ for
all $n>0$. Hence, 
  $$\log \bW(\emptyset,\emptyset,\emptyset) = \frac{(s_1+s_2)(s_1+s_3)(s_2+s_3)}
{s_1s_2s_3} \overline{\bF}_0(q).$$
The equivariant vertex measure is evaluated in 
the Calabi-Yau specialization
in \cite{MNOP1},
$$
\log \bW(\emptyset,\emptyset,\emptyset)|_{s_1+s_2+s_2=0} = M(-q).$$
Hence, $\overline{\bF}_0= -M(-q)$.
\epf

\begin{lm}
\label{divv}
The $q^n$ coefficient of $\bW(\emptyset,\emptyset, \emptyset)$ is
divisible by the cubic factor $(s_1+s_2)(s_1+s_3)(s_2+s_3)$ for all $n>0$.
\end{lm}

\bpf 
We will show the factor $s_1+s_2$ occurs with positive multiplicity in the equivariant vertex 
measure $\bw(\pi)$ for any finite plane partition
$\pi$.  By symmetry, the cyclic permutations of $s_1+s_2$ also occur in $\bw(\pi)$
with positive multiplicity.

Following the
notation of \cite{MNOP1}, let $Q_\pi(t_1,t_2,t_3)$ be the characteristic polynomial of 
the partition $\pi$. Then, the character of the virtual tangent space at $\pi$ is given by
$$\mathrm{V}_{\pi}(t_1,t_2,t_3) = 
Q_\pi - \frac{\bar{Q}_\pi}{t_1t_2t_3} + Q_\pi\bar{Q}_\pi
\frac{(1-t_1)(1-t_2)(1-t_3)}{t_1t_2t_3},$$ 
where $\bar{Q}_\pi(t_1,t_2,t_3) = Q_\pi(t_1^{-1},t_2^{-1},t_3^{-1}).$
The vertex measure is obtained from $\mathrm{V}_{\pi}$ via the prescription 
$$\sum \pm t_1^{i}t_2^{j}t_3^{k} \rightarrow \prod (is_{1}+js_{2}+ks_{3})^{\mp 1}.$$
Hence, the monomials of the form
$t_1^{i}t_2^{i}t_3^{0}$ in $\mathrm{V}_{\pi}$ are 
those which contribute a factor of  $s_{1}+s_{2}$ to $\bw(\pi)$. The total multiplicity 
of $s_1+s_2$ is the negative of 
the constant term in the Laurent polynomial $\mathrm{V}_{\pi}(x, x^{-1},t_3)$.

Let $\rho$ be a 2-dimensional partition. The {\em content} 
of the box $(r,s)$ in $\rho$ is
$r-s$.
The slices of $\pi$ perpendicular to the $z$ direction determine
2-dimensional partitions 
$$\pi_{0},\pi_{1},\pi_2, \ldots\,\  .$$ 
Let $a_{i,j}$ be the number of boxes in $\pi_{j}$ with content $i$.
For convenience, we set $a_{i,j} = 0$ for $j < 0$. 
By definition, 
$Q_\pi(x,x^{-1},t_3) = \sum_{i,j}a_{i,j}x^{i}t_3^{j}$.

The constant term 
of
$\mathrm{V}_{\pi}(x,x^{-1},t_3)$ may be expressed in terms of the contents.
Using
$$ \mathrm{V}_{\pi}(x,x^{-1},t_3)=
Q_\pi(x,x^{-1},t_3) -\frac{\bar{Q}_\pi(x,x^{-1},t_3)}{t_3} 
+Q_\pi\bar{Q}_\pi(2 -x -\frac{1}{x})(\frac{1}{t_3} -1),$$
we find the constant term equals
$$a_{0,0} + 
\sum_{i,j \in {\Z}} (2a_{i,j+1}a_{i,j} - a_{i,j+1}a_{i+1,j}-a_{i+1,j+1}a_{i,j}) 
- (2a_{i,j}a_{i,j}-a_{i,j}a_{i+1,j}-a_{i+1,j}a_{i,j}).$$
We rewrite the constant term in a factored form,
$$ a_{0,0} + \sum_{i,j\in \Z}
\Big( (a_{i,j+1} - a_{i+1,j+1})(a_{i,j}-a_{i+1,j}) - (a_{i,j}-a_{i+1,j})^{2}
\Big)$$
which equals
\begin{equation}
\label{fffq}
a_{0,0} - \frac{1}{2}\sum_{i,j\in \Z}\Big((a_{i,j}-a_{i+1,j}) - (a_{i,j+1}-a_{i+1,j+1})\Big)^{2}.
\end{equation}
Since $(a_{i,0} - a_{i+1,0})= \pm 1$ or $0$, we see 
$$a_{0,0} = \sum_{i \geq 0} (a_{i,0}- a_{i+1,0}) = \sum_{i \geq 0 }(a_{i,0}- a_{i+1,0})^{2}$$
with a similar equality for $i < 0$.  
Therefore, $a_{0,0}$ precisely cancels the $j = -1$ term in \eqref{fffq}, 
yielding the expression
\begin{equation}
\label{rrtt}
-\frac{1}{2}
\sum_{i \in \mathbf{Z}, j \geq 0}((a_{i,j}-a_{i+1,j}) - (a_{i,j+1}-a_{i+1,j+1}))^{2} 
\end{equation}
for the constant term of $\mathrm{V}_{\pi}(x,x^{-1},t_3)$.

We  conclude \eqref{rrtt}
is negative since $a_{i,j} = 0$ for $j$ sufficiently large.
Hence,  the multiplicity of $s_1+s_2$
in $\bw(\pi)$ is strictly positive.
\epf

\begin{Corollary}
The degree 0 localization contributions
over $\infty$ are:
$$\bW_\infty=M(-q)^{\frac{s_2+s_3}{s_1}}.$$
\end{Corollary}

\bpf By Lemmas \ref{gwx} and \ref{llggf}, the Corollary is obtained by extracting
the pole in $s_1$ of $\log \bW(\emptyset,\emptyset,\emptyset)$.
\epf

\subsection{Degree 0 results for toric 3-folds}
Let $X$ be a nonsingular, projective, toric 3-fold equipped with a $\T$-action, and
let $S\subset X$ be a nonsingular toric divisor.

\begin{tm}
$Z_{DT}(X;q)_0 = M(-q)^{\int_X c_3(T_X \otimes K_X)}.$
\end{tm}
\bpf
Let $\{X_\alpha\}$ denote the set of $\T$-fixed points of $X$.
By localization,
$$Z_{DT}(X;q)_0= \prod_{X_\alpha} \ \bW(\emptyset,\emptyset,\emptyset)|_{s_1=-s_1^\alpha,\
s_2=-s_2^\alpha,\ s_3=-s_3^\alpha},$$
where $s_1^\alpha, s_2^\alpha, s_3^\alpha$ are the tangent weights at
$X_\alpha$.
By Theorem \ref{dz},
$$\log Z_{DT}(X;q)_0= \Big( \sum_{X_\alpha} {\frac{(-s^\alpha_1-s^\alpha_2)
(-s^\alpha_1-s^\alpha_3)(-s^\alpha_2-s^\alpha_3)}
{s^\alpha_1s^\alpha_2s^\alpha_3}} \Big)
\cdot \log M(-q).$$
The prefactor on the right is equal to $\int_X c_3(T_X \otimes K_X)$
by a direct application of the Bott residue formula.
\epf

\begin{tm}
$Z_{DT}(X/S;q)_0 = M(-q)^{\int_X c_3(T_X[-S] \otimes K_X[S])}.$
\end{tm}
\bpf
Let $\{S_\gamma\}$ denote the set of $\T$-fixed points of $S$. Let
$s_1^\gamma$ be the normal weight to $S$ at $S_\gamma$, and let 
$s_2^\gamma, s^3_\gamma$ be the tangent weights to $S$ at $S_\gamma$.
By localization,
\begin{eqnarray*}
Z_{DT}(X/S;q)_0 & = & \prod_{X_\alpha\notin S}\ \bW(\emptyset,\emptyset,\emptyset)|_{s_1=-s_1^\alpha,\ 
s_2=-s_2^\alpha,\ s_3=-s_3^\alpha} \ \cdot\\
& & \ \prod_{S_\gamma} \ \bW_\infty|_{s_1=s_1^\gamma,\ s_2=-s_2^\gamma,\ s_3=-s_3^\gamma}\, .
\end{eqnarray*}
By Theorem \ref{dz} and Corollary \ref{dz},
$$\frac{\log Z_{DT}(X/S;q)_0}{\log M(-q)}= 
\sum_{ X_\alpha\notin S}  {\frac{(-s^\alpha_1-s^\alpha_2)
(-s^\alpha_1-s^\alpha_3)(-s^\alpha_2-s^\alpha_3)}
{s^\alpha_1s^\alpha_2s^\alpha_3}} + 
\sum_{S\gamma} {\frac{-s^\gamma_2-s^\gamma_3}{s^\gamma_1}}.$$
The weights of $T_X[-S]\otimes K_X[S] $ at $S_\gamma$ are 
$$-s_2^\gamma-s_3^\gamma,-s^\gamma_2,-s^\gamma_3.$$
Hence, the right side is equal to $\int_X c_3(T_X[-S] \otimes K_X[S])$
by the Bott residue formula.
\epf

\subsection{Evaluations in higher degrees}
While the equivariant vertex measure has a simple formula in degree 1,
$$\bW(1,\emptyset,\emptyset)=   (1+q)^{\frac{s_2+s_3}{s_1}} 
M(-q)^{- \frac{(s_1+s_2)(s_1+s_3)(s_2+s_3)}{s_1s_2s_3}},$$
the higher degree cases are more subtle. We will study 
the evaluations in degrees 1 and higher in a future paper.

%\pagebreak

\vspace{+10 pt}
\noindent
Department of Mathematics \\
Princeton University \\
Princeton, NJ 08544, USA\\
dmaulik@math.princeton.edu \\

\vspace{+10 pt}
\noindent
Institut des Hautes Etudes Scientifiques \\
Bures-sur-Yvette, F-91440, France\\
nikita@ihes.fr \\

\vspace{+10 pt}
\noindent
Department of Mathematics \\
Princeton University \\
Princeton, NJ 08544, USA\\
okounkov@math.princeton.edu \\

\vspace{+10 pt}
\noindent
Department of Mathematics\\
Princeton University\\
Princeton, NJ 08544, USA\\
rahulp@math.princeton.edu

\end{document}